\documentclass{amsart}
\usepackage{amssymb,latexsym,graphics,amsfonts}

\swapnumbers
 \theoremstyle{plain}

 \theoremstyle{definition}
 
 \theoremstyle{remark}
 
 \newcommand{\cal}[1]{\mathcal{#1}}

\begin{document}
\title{Connes amenability of the second dual of Arens regular Banach algebras}
\author{M. Eshaghi Gordji}
\address{Department of Mathematics,
University of Semnan, Semnan, Iran} \email{maj\_ess@Yahoo.com and
meshaghi@semnan.ac.ir}

 \keywords{Derivation , Connes Amenable}
\subjclass[2000]{46HXX}
\dedicatory{}
\smallskip
\begin{abstract}
In this paper we study the Connes amenability of the second dual of
Arens regular Banach algebras. Of course we provide a partial answer
to the question posed by Volker Runde. Also we fined the necessary
and sufficient conditions for the second dual of an Arens regular
module extension Banach algebra to be Connes amenable when the
module is reflexive.
\end{abstract}
\maketitle
$$\bf Introduction$$
A Banach algebra $\cal A$ is said to be dual if there is a closed
submodule ${\cal A}_*$ of $\cal A^*$ such that $\cal A={{\cal
A}_*}^*.$ Let $\cal A$ be a dual Banach algebra. A dual Banach $\cal
A$-module $X$ is called normal if, for every $x\in X,$ the maps
$a\longmapsto a.x$ and
$a\longmapsto x.a$ are $weak^*-weak^*$-continuous from $\cal A$ into $X$. \\
For example if $G$ is a locally compact topological group, then
$M(G)$ is a dual Banach algebra with predual $C_0(G)$. Also if $\cal
A$ is an Arens regular Banach algebra, then $\cal A^{**}$ (by the
first Arens product) is a dual Banach algebra with predual $\cal
A^*.$ Let $\cal A$ be a Banach algebra, and let  $X$ be a Banach
$\cal A$-module then a derivation from $\cal A$ into $X$ is a linear
map $D$, such that for every $a,b \in \cal A,$
$D(ab)=D(a).b+a.D(b).$ Let  $x\in X,$ and let $\delta_x:\cal
A\longrightarrow X$ defined by $\delta_x(a)=a.x-x.a\hspace
{0.5cm}(a\in \cal A),$ then $\delta_x$ is a derivation, derivations
of this form are called inner derivations. A Banach algebra is
called amenable if every derivation from $\cal A$ into each dual
$\cal A$-module is inner; i.e. $H^1(\cal A,X^*)=\{o\}$, foe every
$\cal A$-module $X$. This definition was introduced by B. E. Johnson
in [4]. A dual Banach algebra $\cal A$ is Connes amenable if every
$weak^*-weak^*$-continuous derivation from $\cal A$ into each normal
dual Banach $\cal A$-module $X$
 is inner; i.e. $H^1_{w^*}(\cal A, X)=\{o\}$, this definition was introduced by V.Runde
(see section 4 of [6]).
We answer partially to the following question [6, 2]. \\
Let $\cal A$ be an Arens regular Banach algebra such that $\cal
A^{**}$ is Connes amenable. Need $\cal A$ be amenable?
 \section{Second dual of Arens regular Banach algebras}
Let the  Banach algebra $\cal A$ be Arens regular, then we have the following assertions.\\
(i) If $\cal A$ is  amenable, then  $\cal A^{**}$ is Connes
amenable.\\
 (ii) If $\cal A$ is an
ideal of $\cal A^{**}$ and  $\cal A^{**}$ is Connes amenable, then
$\cal A$ is amenable [5].\\
First we have the following result.
\paragraph{\bf Theorem 1.1.}
Let $\cal A$ be a Banach algebra which $A^{**}$ is Arens regular
and $\cal A^{****}$ is Connes amenable. Then $\cal A^{**}$ is
Connes amenable. Also if $\cal A^{**}$ is an ideal of $\cal
A^{****}$, then $\cal A^{**}$ is amenable and reflexive.
\paragraph{\bf Proof.}
Let $X$ be a normal $\cal A^{**}$-module, and let
$\pi:a''''\longmapsto a''''\mid_{\cal A^{*}}: \cal A^{****}
\longrightarrow \cal A^{**}$ be the restriction  map. Since $\pi$
is $weak^*-weak^*$-continuous, then $X$ is a normal $\cal
A^{****}$-module by the following module actions
$$a''''x=\pi(a'''')x, \hspace{0.7cm} xa''''=x\pi(a'''')\hspace{0.7cm} (x\in X, a''''\in \cal A^{****}).$$
Let $D:\cal A^{**}\longrightarrow X$ be a $weak^*-weak^*$-continuous
derivation. It is easy to show that $Do\pi:\cal
A^{****}\longrightarrow X$ is a $weak^*-weak^*$-continuous
derivation. If $\cal A^{****}$ is Connes amenable then $Do\pi$ is
inner, so $D$ is inner, and $\cal A^{**}$ is Connes amenable. Connes
amenability of $\cal A^{**}$ implies that $\cal A^{**}$ is unital.
Let now $\cal A^{**}$ be an ideal of $\cal A^{****}$, then $\cal
A^{**}=\cal A^{****}$. Thus $\cal A^{**}$ is reflexive, also by
theorem 4.4 of [5], $\cal A^{**}$ is amenable. \hfill$\blacksquare~$

\paragraph{\bf Corollary 1.2.}
Let $\cal A$ be a Banach algebra which $A^{**}$ is amenable and
Arens regular. If $\cal A^{**}$ is an ideal of $\cal A^{****}$,
then $\cal A$ is reflexive.

\paragraph{\bf Theorem 1.3.} Let $\cal A$ be an Arens regular Banach
algebra with a bounded approximate identity, which is a right
ideal of $A^{**}$. Let for every $\cal A^{**}$-neo unital module
$X$, $X^*$ factors $\cal A$ on the left, i.e. $\cal AX^*=X^*$. If
$\cal A^{**}$ is Connes amenable, then $\cal A$ is amenable.
\paragraph{\bf Proof.}
Let $X$ be a $\cal A$-module, and let $D:\cal A\longrightarrow
X^{**}$ be a derivation. Since  $\cal A^{**}$ is Connes amenable,
then $\cal A^{**}$ has unite element $E$. We can extend the actions
of ${\cal A}$ on $X^{**}$ to actions of ${\cal A}^{**}$ on
$X^{****}$ via
\[a''.x'''' =w^{*}\textup{-}\lim_{i} \lim_{j} a_{i} \, x''_{j}
\]
and
\[ x''''.a'' =w^{*}\textup{-}\lim_{j} \lim_{i} x''_{j} \, a_{i},\]
where $a'' =w^{*}\textup{-}\lim_{i} a_{i}$, \  $x''''
=w^{*}\textup{-}\lim_{j} x''_{j}$. We have the direct sum
decomposition $$X^{****}=EX^{****}E \oplus (1-E)X^{****}E \oplus
EX^{****}(1-E) \oplus (1-E)X^{****}(1-E),$$ as $\cal
A^{**}$-modules. For i=1,2,3,4, let $\pi_i$ be the associated
projection and let $D_i=\pi_ioD^{**}$. $\pi_i$ is a $\cal
A^{**}$-module morphism, then $D_i$ is a derivation. Let
$y_2=-D_2(E)$. Since $D^{**}(a'')=(D^{**}(E))a''-ED^{**}(a'')$ and
$a''D_2(E)=a''(1-E)(D^{**}(E))E=0$,  then
$$D_2(a'')=(1-E)(D^{**}(E))a''E=D_2(E)a''=\delta_{y_2}(a'').$$
A similar argument applies to $D_3$ and $D_4$.  We show that
$EX^{****}E$ is a normal $\cal A^{**}$-module. First we have
$$EX^{****}E=(EX^{***}E)^* \hspace {0.7cm} (1).$$ Let $Ex''''E\in EX^{****}E$ and let
$a''_\alpha~\stackrel{weak^*}{-\hspace{-.2cm}-\hspace{-.2cm}\longrightarrow}
a''$ in $\cal A^{**}$. Since $EX^{***}E$ is neo-unital $\cal
A^{**}$- module, then $ EX^{****}E=EX^{****}E\cal A$, therefore
there is $a\in \cal A$ and $y''''\in X^{****}$ such that
$Ey''''Ea=Ex''''E.$ We have
$aa''_\alpha~\stackrel{weak^*}{-\hspace{-.2cm}-\hspace{-.2cm}\longrightarrow}
aa''$ in $\cal A^{**},$ since $\cal A$ is a right ideal of
$A^{**}$, then
$aa''_\alpha~\stackrel{weak}{-\hspace{-.2cm}-\hspace{-.2cm}\longrightarrow}
aa''$ in $\cal A.$ Thus by (1), we have $$Ex''''Ea''_\alpha
=Ey''''Eaa''_\alpha~\stackrel{weak}{-\hspace{-.2cm}-\hspace{-.2cm}\longrightarrow}
Ey''''Eaa''=Ex''''Ea'' \hspace{1cm} in \hspace{0.3cm}
EX^{****}E.$$ Then $$Ex''''Ea''_\alpha
~\stackrel{weak^*}{-\hspace{-.2cm}-\hspace{-.2cm}\longrightarrow}
Ex''''Ea'' \hspace{1cm} in \hspace{0.3cm} EX^{****}E.$$ Trivially
we have $$a''_\alpha Ex''''E
~\stackrel{weak^*}{-\hspace{-.2cm}-\hspace{-.2cm}\longrightarrow}
a''Ex''''E \hspace{1cm} in \hspace{0.3cm}  EX^{****}E.$$ This
means that $EX^{****}E$ is a normal $\cal A^{**}$-module. Since
$\pi_1$ and $D^{**}$ are $weak^*-weak^*$-continuous, and $\cal
A^{**}$ is Connes amenable, then there is $x_1''''\in X^{****}$
such that $D_1=\delta_{Ex_1''''E}$. Thus $D^{**}$ is inner. On the
other hand we have the following direct sum decomposition of $\cal
A$-moduls
$$X ^{****}=\widehat{X^{**}}\oplus {\widehat {(X^{*})}}^{\perp}.$$
Let $\pi:X^{****}\longrightarrow X^{**}$ be the natural projection,
then $D=\pi oD^{**}$ is inner. Thus $H^1(\cal A, X^{**})=\{o\}$, and
by Proposition 2.8.59 of [1], $\cal A$ is
amenable.\hfill$\blacksquare~$
\section{Module extension dual Banach algebras}
Let $\cal A$  be a Banach algebra and $M$ be a Banach $\cal
A$-module (with module actions $\pi_r$ and $\pi_l$), let ${\cal
B}=M\oplus_1 {\cal A}$ as a Banach space, so that
$$\|(m,a)\|=\|m\|+\|a\|\hspace{1.5cm} (a\in{\cal A}~,~m\in M)~.$$
Then $\cal B$ is a Banach algebra with the product
$$(m_1,a_1)(m_2,a_2)=(m_1\cdot a_2+a_1\cdot m_2,a_1a_2)~.$$
The second dual ${\cal  B}^{**}$ of $\cal B$ is identified  with
$M^{**}\oplus_1 {\cal A}^{**}$ as a Banach space and the first
Arens product on ${\cal B}^{**}$ is given by
$$(m''_1,a'' _1)(m''_2,a'' _2)=(m''_1\cdot
a'' _2+a'' _1\cdot m''_2,a''_1 a''_2)~.$$
 As in [3] we can show that
$B$ is Arens regular if and only if for every $a''\in \cal
A^{**}$, and $m''\in M^{**},$ \\
(1) $b''\longmapsto a'' b'':{\cal A}^{**}\longrightarrow {\cal
A}^{**}$ is
$weak^*-weak^*$ continuous.\\
(2) $n''\longmapsto a''n'':M^{**}\longrightarrow M^{**}$ is
$weak^*-weak^*$ continuous.\\
(3) $b''\longmapsto m''b'':{\cal A}^{**}\longrightarrow M^{**}$ is
$weak^*-weak^*$ continuous.\\
Then we have the following theorem.

\paragraph{\bf Theorem 2.1.} Let $\cal A$ be
an Arens regular Banach algebra, and let $M$ be a reflexive Banach $\cal A$-module. Then\\
(i) $\cal B=M\oplus_1 {\cal A}$ is Arens regular.\\
(ii) $\cal B^{**}={(M\oplus_1 {\cal A})}^{**}$ is Connes amenable if
and only if $M=0$, and $\cal A^{**}$ is Connes amenable.
\paragraph{\bf Proof.} We can prove (i) by the argument above theorem.
To prove (ii), suppose that $\cal B^{**}$ is Connes amenable, we
need only to show that $M=0$. Let $X=M^{***}\bigotimes_p {\cal
A}^{**}$. We define the module actions of $\cal B$ on $X$ as
follows:
$$(m'''\otimes_pa'').(b'',x'')=m'''\otimes_pa''b'', \hspace{0.7cm}(b'',x'')
.(m'''\otimes_pa'')=b''m'''\otimes_pa'',$$

so we define $D:\cal B^{**}\longrightarrow X^*$ by
$$\langle D((a'',x'')),m'''\otimes b''\rangle=\langle x''m''',b''\rangle.$$
Where $m'''\in M^{***}, x''\in M^{**}$ and $a'',b''\in \cal A^{**}.$

Let
$(b''_\alpha,x''_\alpha)~\stackrel{weak^*}{-\hspace{-.2cm}-\hspace{-.2cm}\longrightarrow}
(b'',x'')$ in $\cal B^{**},$ then we have
$b''x''_\alpha~\stackrel{weak^*}{-\hspace{-.2cm}-\hspace{-.2cm}\longrightarrow}
b''x''$ in $M^{**}.$ Since $M$ is reflexive then
$b''x''_\alpha~\stackrel{weakly}{-\hspace{-.2cm}-\hspace{-.2cm}\longrightarrow}
b''x''$ in $M^{**}.$ Then for every $m'''\in M^{***}$ we have
$\langle x''_\alpha m''',b''\rangle\longrightarrow \langle x''
m''',b''\rangle.$ This means that $D$ is $weak^*-weak^*$ continuous.
Also for every $(a''_1,x''_1) , (a''_2,x''_2) \in \cal B^{**},
m'''\in M^{***}$ and $ x''\in M^{**},$  we have
\begin{align*}
\langle D((a''_1,x''_1)(a''_2x''_2)),m'''\otimes b''\rangle=&
\langle D((a''_1a''_2,x''_1a''_2+a''_1x''_2)),m'''\otimes
b''\rangle\\
=& \langle (a''_1x''_2m'''+x''_1a''_2m'''),b''\rangle=
\langle x''_2m''',b''a''_1\rangle+\langle x''_1a''_2m''',b''\rangle\\
=&D((a''_1,x''_1)),(a''_2x''_2).(m'''\otimes b'')\rangle+
D((a''_2x''_2)),(m'''\otimes b'').(a''_1,x''_1)\rangle\\
=&D((a''_1,x''_1)).(a''_2x''_2),m'''\otimes b''\rangle+
(a''_1,x''_1).D((a''_2x''_2)),m'''\otimes b''\rangle.
\end{align*}
Thus $D$ is a derivation. Since $\cal B^{**}$ is Connes amenable,
than it is unital. Let $(E,x''_1)$ be the unit element of $\cal
B^{**}$, then it is easy to show that $x''_1=0$ and that $E$ is unit
element of $\cal A^{**}$, so $Em''=m''E=m''$ and $Em'''=m'''E=m''',$
for every $m''\in M^{**}$ and $m'''\in M^{***}.$ Since $M$ is
reflexive then it is easy to show that $\cal A M=M\cal A=M$ and the
module actions of $\cal B^{**}$ on $X^*$ are $weak^*-weak^*$
continuous. Then $D$ is inner, and there exists $F\in X^*$ such that
$D(a'',m'')=(a'',m'').F-F.(a'',m'').$ Then for every $m'''\otimes
b''\in X,$ we have

\begin{align*}
\langle m''m''',b''\rangle=& \langle \langle
D((a'',m'')),m'''\otimes b''\rangle\\
=& \langle F,(m'''\otimes b'').(a'',m'')\rangle+\langle
F,(a'',m'').(m'''\otimes
b'')\rangle\\
=&\langle F,m'''\otimes b''a''+a''m'''\otimes b''\rangle.
\end{align*}
Let $a''=0,$ then we have $m''m'''=0$ for every $m'''\in M^{***},
m''\in M^{**}.$ This means that $\cal A^{**}M^{**}=0$. Thus
$m''=Em''=0$ for every $m''\in M^{**}$ and the proof is
complete.\hfill$\blacksquare~$

Let $\cal A$ be a Banach algebra and let $\varphi \in \Omega_{\cal
A}$ be a multiplier on $\cal A$. Then $\Bbb C$ is a Banach $\cal
A$-module by module actions
$$a.c=\varphi(a)c, ~~~~~~~~~~c.a=c\varphi(a),  ~~~~~~~~(a\in \cal A , c\in
\Bbb C).$$ We denote this $\cal A$-module with $\Bbb C_{\varphi}$.
By apply above theorem we can give a class of dual Banach algebras
which are not Conns amenable.
\paragraph{\bf Corollary 2.2.} Let $\cal A$ be
an Arens regular Banach algebra, and let $0\neq \varphi \in
\Omega_{\cal A}$. Then ${(\Bbb C_{\varphi} \oplus_1 {\cal A})}^{**}$
is a dual Banach algebra which is not Connes amenable.


\end{document}